\input amstex
\magnification=\magstep1
\documentstyle{amsppt}
\pagewidth{16 true cm} \pageheight{22.0 true cm}
\TagsOnLeft\CenteredTagsOnSplits
\input epsf.sty

\def\ignore#1{\relax}
\def\f{\varphi}
\def\t{\tau}
\def\yf{\Bbb{YF}}
\def\Q{\Bbb Q}
\def\Fun{\hbox{\rm Fun}}
\def\Y{\Bbb Y}

\topmatter
\title The Martin Boundary of the Young-Fibonacci Lattice \endtitle
\author
{Frederick M. Goodman} {\rm   and  }
{Sergei V. Kerov}
\endauthor
\address Department of Mathematics,
University of Iowa,
Iowa City, Iowa 52242, USA
\endaddress
\address Steklov Math. Institute (POMI),
Fontanka 27,
St.~Petersburg 191011, Russia
\endaddress
\abstract
In this paper we find the Martin boundary for the Young-Fibonacci
lattice $\Bbb{YF}$. Along with the lattice of Young diagrams,
this is the most interesting example of a differential partially
ordered set. The Martin boundary construction provides an
explicit Poisson-type integral representation of non-negative
harmonic functions on $\Bbb{YF}$. The latter are in a canonical
correspondence with a set of traces on the locally semisimple Okada
algebra. The set is known to contain all the indecomposable
traces. Presumably, all of the traces in the set are
indecomposable, though we have no proof of this conjecture. Using
an explicit product formula for Okada characters, we derive
precise regularity conditions under which a sequence of
characters of finite-dimensional Okada algebras converges.
\endabstract
\thanks
The second author was partially supported by the grant
INTAS-94-3420
\endthanks
\keywords
Differential poset, harmonic functions, Martin boundary, Okada
algebras, non-commutative symmetric functions
\endkeywords
\toc
\widestnumber\head{\S 9}
\head \S 1. Introduction \endhead
\head \S 2. The Young -- Fibonacci lattice\endhead
\head \S 3. Harmonic functions on graphs and traces of $AF$-algebras\endhead
\head \S 4. Harmonic functions on differential posets\endhead
\head \S 5. Okada clone of the symmetric function ring\endhead
\head \S 6. A product formula for Okada characters\endhead
\head \S 7. The Martin boundary of the Young -- Fibonacci lattice\endhead
\head \S 8. Regularity conditions\endhead
\head \S 9. Concluding remarks\endhead
\specialhead  Appendix\endspecialhead
\specialhead  References\endspecialhead
\endtoc
\endtopmatter

\document

\head \S 1. Introduction \endhead

The Young-Fibonacci lattice $\yf$ is a fundamental example of a
{\it differential partially ordered set} which was introduced by
R.~Stanley \cite{St1} and S.~Fomin [F1]. In many ways, it is
similar to another major example of a differential poset, the
Young lattice $\Bbb{Y}$. Addressing a question posed by Stanley,
S.~Okada has introduced \cite{Ok} two algebras associated to
$\yf$. The first algebra $\Cal F$ is a locally semisimple algebra
defined by generators and relations, which bears the same
relation to the lattice $\yf$ as does the group algebra
$\Bbb{C}\frak{S}_\infty$ of the infinite symmetric group to
Young's lattice. The second algebra $R$ is an algebra of
non-commutative polynomials, which bears the same relation to the
lattice $\yf$ as does the ring of symmetric functions to Young's
lattice.

The purpose of the present paper is to study some combinatorics,
both finite and asymptotic, of the lattice $\yf$.  Our object of
study is the compact convex set of {\it harmonic functions} on
$\yf$ (or equivalently the set of positive normalized traces on
$\Cal F$ or certain positive  linear functionals on $R$.)  We
address the study of harmonic functions by determining the {\it
Martin boundary} of the lattice $\yf$. The Martin boundary is the
(compact)  set consisting of those harmonic functions which can
be obtained by finite rank approximation. There are two basic
facts related to the Martin boundary construction:\quad 1) every
harmonic function is represented by the integral of a probability
measure on the Martin boundary, and\quad 2) the set of extreme
harmonic functions is a subset of the Martin boundary (see, e.g.,
\cite{D}).

This paper gives a parametrization of the Martin boundary for
$\yf$  and a description of its topology.

The Young-Fibonacci lattice is described in Section 2, and
preliminaries on harmonic functions are explained in Section 3.
A first rough description of our main results is given at the
end of Section 3. (A precise description of the parametrization
of harmonic functions is  found in Section 7, and the proof,
finally, is contained in Section 8.) Section 4 contains some
general results on harmonic functions on differential posets.

The main tool in our study is the Okada ring $R$ and two bases of
this ring, introduced by Okada, which are in some respect
analogous to  the Schur function basis and the power sum function
basis in the ring of symmetric functions (Section 5). We describe
the Okada analogs of the Schur function basis by non-commutative
determinants of tridiagonal matrices with monomial entries. We
obtain a simple and explicit formula for the transition matrix
(character matrix) connecting  the s-basis and the p-basis, and
also for the value of (the linear extension of) harmonic
functions evaluated on the p-basis. This is done in Sections 6
and 7.

The explicit formula allows us to study the regularity question
for the lattice $\yf$, that is the question of convergence of
extreme traces of  finite dimensional Okada  algebras $\Cal F_n$
to traces of the inductive limit algebra $\Cal F=\varinjlim\Cal
F_n$.  The regularity question is studied in Section 8.

The analogous questions for Young's lattice $\Bbb Y$ (which is
also a differential poset) were answered some time ago. The
parametrization of the Martin boundary of $\Bbb Y$ has been
studied in \cite{Th}, \cite{VK}. A different approach was
recently given in \cite{Oku}.

A remaining open problem for the Young-Fibonacci lattice is to
characterize the set of extreme harmonic functions within the
Martin boundary. For Young's lattice, the set of extreme harmonic
functions coincides with the entire Martin boundary.

\smallskip\noindent
{\bf Acknowledgement.} The second author thanks the Department of
Mathematics, University of Iowa, for a teaching position in the
Spring term of 1993, during which  most  of this work was done.
This paper was completed in May, 1997 at the home of Sergey and Irina
Fomin, whom we thank for their most generous hospitality.

\head \S 2. The Young-Fibonacci lattice \endhead

In this Section we recall the definition of Young-Fibonacci
modular lattice (see Figure 1) and some basic facts related to
its combinatorics. See Section A.1 in the Appendix for the
background definitions and notations related to graded graphs and
differential posets. We refer to \cite{F1-2}, \cite{St1-3} for a
more detailed exposition.

\subhead  A simple recurrent construction \endsubhead

The simplest way to define the graded graph
$\Bbb{YF}=\bigcup_{n=0}^\infty\Bbb{YF}_n$ is provided by the
following recurrent procedure.

Let the first two levels $\Bbb{YF}_0$ and $\Bbb{YF}_1$ have just
one vertex each, joined by an edge. Assuming that the part of the
graph $\Bbb{YF}$, up to the $n$th level $\Bbb{YF}_n$, is already
constructed, we define the set of vertices of the next level
$\Bbb{YF}_{n+1}$, along with the set of adjacent edges, by first
{\it reflecting} the edges in between the two previous levels,
and then by {\it attaching} just one new edge leading from each
of the vertices on the level $\Bbb{YF}_n$ to a corresponding new
vertex at level $n+1$.

$$
\epsffile{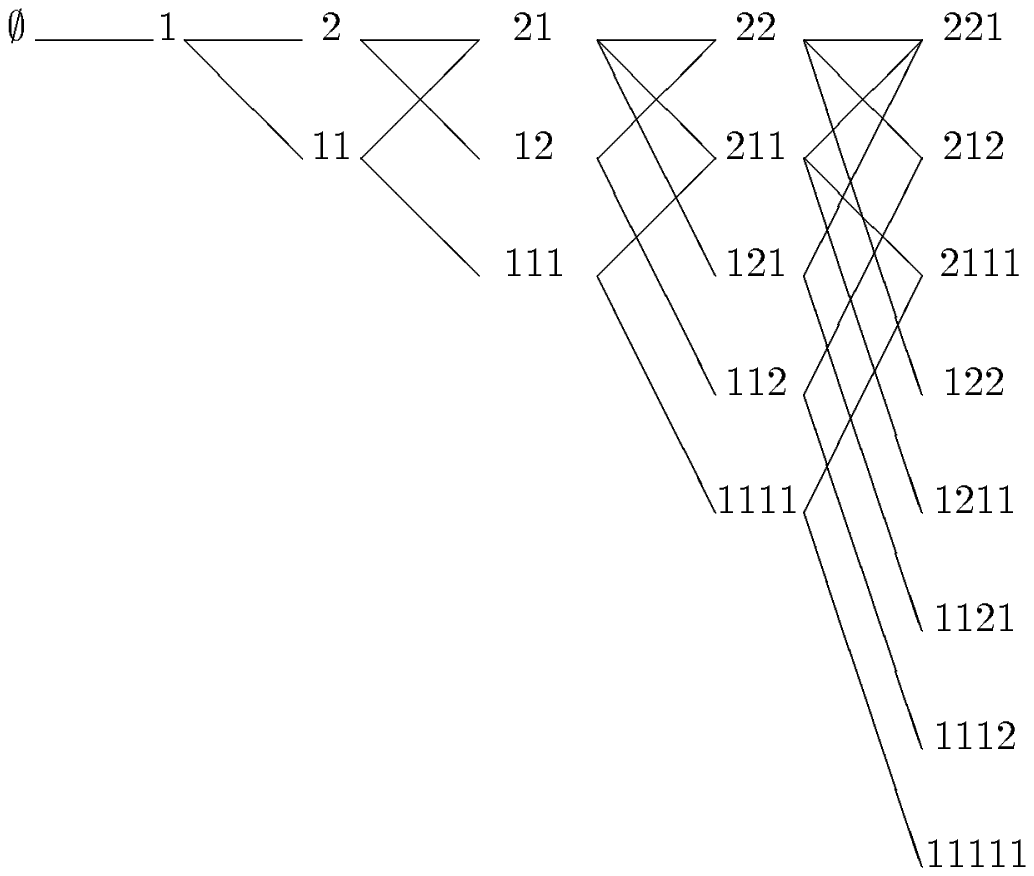}
$$
\centerline{Figure 1.\quad The Young-Fibonacci lattice.}
\bigskip

In particular, we get two vertices in the set $\Bbb{YF}_2$, and
two new edges: one is obtained by reflecting the only existing
edge, and the other by attaching a new one. More generally, there
is a natural notation for new vertices which helps to keep track
of the inductive procedure. Let us denote the vertices of
$\Bbb{YF}_0$ and $\Bbb{YF}_1$ by an empty word $\varnothing$ and
$1$ correspondingly. Then the endpoint of the reflected edge will
be denoted by $2$, and the end vertex of the new edge by $11$. In
a similar way, all the vertices can be labeled by words in the
letters $1$ and $2$. If the left (closer to the root
$\varnothing$) end of an edge is labeled by a word $v$, then the
endvertex of the reflected edge is labeled by the word $2v$. Each
vertex $w$ of the $n$th level is joined to a vertex $1w$ at the
next level by a new edge (which is not a reflection of any
previous edge).

Clearly, the number of vertices at the $n$th level $\Bbb{YF}_n$
is the $n$th Fibonacci number $f_n$.

\subhead Basic definitions \endsubhead

We now give somewhat more formal description of the
Young-Fibonacci lattice and its Hasse diagram.

\smallskip\noindent
{\bf Definition.}
A finite word in the two-letter alphabet $\{1,2\}$ will be
referred to as a {\it Fibonacci word}. We denote the sum of
digits of a Fibonacci word $w$ by $|w|$, and we call it the {\it
rank} of $w$. The set of words of a given rank $n$ will be
denoted by $\Bbb{YF}_n$, and the set of all Fibonacci words by
$\Bbb{YF}$. The {\it head} of a Fibonacci word is defined as the
longest contiguous subword of $2$'s at its left end. The {\it
position} of a 2 in a Fibonacci word is one more than the rank of
the subword to the right of the 2; that is if $w = u2v$, then the
position of the indicated 2 is $|v| + 1$.

Next we define a partial order on the set $\Bbb{YF}$ which is
known to make $\Bbb{YF}$ a modular lattice. The order will be
described by giving  the covering relations on $\Bbb{YF}$ in two
equivalent forms.

Given a Fibonacci word $v$, we first define the set
$\overline{v}\subset\Bbb{YF}$ of its successors. By definition,
this is exactly the set of words $w\in\Bbb{YF}$ which can be
obtained from $v$ by one of the following three operations:

(i)\hskip4mm put an extra $1$ at the left end of the word $v$;

(ii)\hskip3mm replace the first $1$ in the word $v$ (reading left
to right) by $2$;

(iii)\hskip2mm insert $1$ anywhere in between $2$'s in the head
of the word $v$, or immediately after the last $2$ in the head.

\noindent {\bf Example.}
Take $222121112$ for the word $v$ of rank $14$. Then the group of
3 leftmost $2$'s forms its head, and $v$ has 5 successors, namely
$$
\overline{v}=
\{\bold{1}222121112,\,
2\bold{1}22121112,\,
22\bold{1}2121112,\,
222\bold{1}121112,\,
222\bold{2}21112\}.
$$
The changing letter is shown in boldface. Note that the ranks of
all successors of a Fibonacci word $v$ are one bigger than that
of $v$.

The set $\underline{v}$ of predecessors of a non-empty Fibonacci
word $v$ can be described in a similar way. The operations to be
applied to $v$ in order to obtain one of its predecessors are as
follows:

(i)\hskip 4mm the leftmost letter $1$ in the word $v$ can be
removed;

(ii)\hskip 3mm any one of $2$'s in the head of $v$ can be
replaced by $1$.

\noindent {\bf Example.}
The word $v=222121112$ has 4 predecessors, namely
$$
\underline{v}=\{
122121112,\,
212121112,\,
221121112,\,
22221112\}.
$$

We write $u\nearrow v$ to show that $v$ is a successor of $u$
(and $u$ is a predecessor of $v$). This is a covering relation
which determines a partial order on the set $\Bbb{YF}$ of
Fibonacci words. As a matter of fact, it is a modular lattice,
see \cite{St1}. The initial part of the Hasse diagram of the
poset $\Bbb{YF}$ is represented in Figure 1.

\subhead The Young-Fibonacci lattice as a differential poset
\endsubhead

Assuming that the head length of $v$ is $k$, the word $v$ has
$k+2$ successors and $k+1$ predecessors, if $v$ contains at least
one letter $1$. If $v=22\ldots2$ is made of $2$'s only, it has
$k+1$ successors and $k$ predecessors. Note that the number of
successors is always one bigger than that of predecessors.
Another important property of the lattice $\Bbb{YF}$ is that, for
any two different Fibonacci words $v_1$, $v_2$ of the same rank,
the number of their common successors equals that of common
predecessors (both numbers can only be $0$ or $1$). These are
exactly the two characteristic properties (D1), (D2) of
differential posets, see Section A.1. In what follows we shall
frequently use the basic facts on differential posets, surveyed
for the reader's convenience in the Appendix. Much more
information on differential posets and their generalizations can
be found in \cite{F1}, \cite{St1}.

\subhead The Okada algebra \endsubhead

Okada [Ok] introduced a (complex locally semisimple) algebra
$\Cal F$, defined by generators and relations, which admits the
Young-Fibonacci lattice $\Bbb{YF}$ as its branching diagram. The
Okada algebra has generators $(e_i)_{i \ge 1}$ satisfying the
relations:

\medskip\noindent
$({\text O1})\hskip 1.5em e_i^2=e_i$  \quad for all $i\ge1$;
\smallskip\noindent
$({\text O}2)\hskip .9em e_{i}e_{i-1}e_{i}=  {1\over i}\;e_{i}$
\quad for all $i\ge2$;
\smallskip\noindent
$({\text O}3)\hskip 1em e_ie_j=e_je_i$ \quad for $|i-j|\geq 2$.
\smallskip\noindent

The algebra $\Cal F_n$ generated by the first $n-1$ generators
$e_1, \dots, e_{n-1}$ and these identities is semisimple of
dimension $n!$, and has simple modules $M_v$ labelled by elements
$v \in \Bbb{YF}_n$. For $u \in \yf_{n-1}$ and $v \in \yf_n$, one
has $u \nearrow v$ if, and only if, the simple $\Cal F_n$-module
$M_v$, restricted to the algebra $\Cal F_{n-1}$ contains the
simple $\Cal F_{n-1}$-module $M_u$. As a matter of fact, the
restrictions of simple $\Cal F_n$-modules to $\Cal F_{n-1}$ are
multiplicity free.

\head \S 3. Harmonic functions on graphs and traces of
$AF$-algebras \endhead

In this Section,   we recall the notion of harmonic functions on
a {\it graded graph} and the classical Martin boundary
construction for graded graphs and {\it branching diagrams}. We
discuss the connection between harmonic functions on branching
diagrams and traces on the corresponding $AF$-algebra.  Finally,
we give a preliminary statement on our main results on the Martin
boundary of the Young-Fibonacci lattice.

We refer the reader to Appendix A.1 for basic definitions on
graded graphs and branching diagrams and to \cite{E}, \cite{KV}
for more details on the combinatorial theory of $AF$-algebras.

\subhead The Martin boundary of a graded graph \endsubhead

A function $\varphi:\Gamma\to\Bbb{R}$ defined on the set of
vertices of a graded graph $\Gamma$ is called {\it harmonic}, if
the following variant of the ``mean value theorem'' holds for all
vertices $u \in\Gamma$:
$$
\varphi(u ) = \sum_{w :u \nearrow w } \varphi(w ).
\tag 3.1
$$
We are interested in the problem of determining the space
$\Cal{H}$ of all non-negative harmonic functions normalized at
the vertex $\varnothing$ by the condition
$\varphi(\varnothing)=1$. Since $\Cal{H}$ is a compact convex set
with the topology of pointwise convergence, it is interesting to
ask about its set of extreme points.

A general approach to the problem of determining the set of
extreme points is based on the Martin boundary construction (see,
for instance, \cite{D}). One starts with the {\it dimension
function} $d(v ,w )$ defined as the number of all oriented paths
from $v $ to $w $. We put $d(w )=d(\varnothing,w )$.

From the point of view of potential theory, $d(v ,w )$ is the
Green function with respect to ``Laplace operator''
$$
(\Delta\varphi)(u ) = -\varphi(u ) +
\sum_{w :u \nearrow w } \varphi(w ).
\tag 3.2
$$
This means that if $\psi_w (v )=d(v ,w )$ for a fixed vertex $w
$, then $-(\Delta\psi_w )(v )=\delta_{v w }$ for all $v
\in\Gamma$. The ratio
$$
K(v ,w ) = \frac{d(v ,w )}{d(w )}
\tag 3.3
$$
is usually called the {\it Martin kernel}.

Consider the space $\Fun(\Gamma)$ of all functions
$f:\Gamma\to\Bbb{R}$ with the topology of pointwise convergence,
and let $\widetilde{E}$ be the closure of the subset
$\widetilde{\Gamma}\subset\Fun(\Gamma)$ of functions $v
\mapsto{K}(v ,w )$, $w \in\Gamma$. Since those functions are
uniformly bounded, $0\le{K}(v ,w )\le1$, the space
$\widetilde{E}$ (called the {\it Martin compactification}) is
indeed compact. One can easily check that
$\widetilde{\Gamma}\subset\widetilde{E}$ is a dense open subset
of $\widetilde{E}$. Its boundary
$E=\widetilde{E}\setminus\widetilde{\Gamma}$ is called the {\it
Martin boundary} of the graph $\Gamma$.

By definition, the Martin kernel (3.3) may be extended by
continuity to a function $K:\Gamma\times\widetilde{E}\to\Bbb{R}$.
For each boundary point $\omega\in{E}$ the function
$\varphi_\omega(v )=K(v ,\omega)$ is non-negative, harmonic, and
normalized. Moreover, harmonic functions have an integral
representation similar to the classical Poisson integral
representation for non-negative harmonic functions in the disk:

\proclaim{Theorem \rm (cf. \cite{D})}
Every normalized non-negative harmonic function
$\varphi\in\Cal{H}$ admits an integral representation
$$
\varphi(u ) = \int_E K(u ,\omega)\;M(d\omega),
\tag 3.4
$$
where $M$ is a probability measure. Conversely, for every
probability measure $M$ on $E$, the integral (3.4) provides a
non-negative harmonic function $\varphi\in\Cal{H}$.
\endproclaim

All {\it indecomposable} (i.e., {\it extreme}) functions in
$\Cal{H}$ can be represented in the form $\varphi_\omega(v )=K(v
,\omega)$, for appropriate boundary point $\omega\in{E}$, and we
denote by $E_{min}$ the corresponding subset of the boundary $E$.
It is known that $E_{min}$ is a non-empty $G_\delta$ subset of
$E$. One can always choose the measure $M$ in the integral
representation (3.4) to be supported by $E_{min}$. Under this
assumption, the measure $M$ representing a function
$\varphi\in\Cal{H}$ via (3.4) is unique.

Given a concrete example of a graded graph, one looks for an
appropriate ``geometric'' description of the abstract Martin
boundary. The purpose of the present paper is to give an explicit
description for the Martin boundary of   the Young -- Fibonacci
graph $\Bbb{YF}$.

\subhead The traces on locally semisimple algebras \endsubhead

We next discuss the relation between harmonic functions on a
graded graph and traces on locally semisimple algebras.  A {\it
locally semisimple complex algebra $A$} (or AF-algebra) is the
union of an increasing sequence of finite dimensional semisimple
complex algebras, $A=\varinjlim A_n$. The {\it branching diagram}
or {\it Bratteli diagram} $\Gamma(A)$ of a locally semisimple
algebra $A$ (more precisely, of the approximating sequence
$\{A_n\}$) is a graded graph whose vertices of rank $n$
correspond to the simple $A_n$-modules. Let $M_v$ denote the
simple $A_n$-module corresponding to a vertex $v\in\Gamma_n$.
Then a vertex $v$ of rank $n$ and a vertex $w$ of rank $n+1$ are
joined by $\varkappa(v,w)$ edges if the simple $A_{n+1}$ module
$M_w$, regarded as an $A_n$ module, contains $M_v$ with
multiplicity $\varkappa(v,w)$.  We will assume here that all
multiplicities $\varkappa(v,w)$ are 0 or 1, as this is the case
in the example of the Young-Fibonacci lattice with which we are
chiefly concerned. Conversely, given a branching diagram $\Gamma$
-- that is, a graded graph with unique minimal vertex at rank 0
and no maximal vertices -- there is a locally semisimple algebra
$A$ such that $\Gamma(A) = \Gamma$.

A {\it trace} on a locally semisimple algebra $A$ is a complex
linear functional $\psi$ satisfying
$$
\aligned
\psi(e) &\ge 0 \quad \text{for all idempotents}\quad e\in A; \\
\psi(1) &=1; \\
\psi(ab) &= \psi(ba) \quad \text{for all} \quad a,b \in A.
\endaligned
\tag 3.5
$$
To each trace $\psi$ on $A$, there corresponds a positive
normalized harmonic function $\tilde \psi$ on $\Gamma =
\Gamma(A)$ given by
$$
\tilde \psi(v) = \psi(e)
\tag 3.6
$$
whenever $v$ has rank $n$ and $e$ is a minimal idempotent in
$A_n$ such that $e M_v \ne (0)$ and $e M_w = 0$ for all $w \in
\Gamma_n \setminus \{v\}$. The trace property of $\psi$ implies
that $\tilde \psi$ is a well defined non-negative function on
$\Gamma$, and harmonicity of $\tilde \psi$ follows from the
definition of the branching diagram  $\Gamma(A)$.  Conversely, a
positive normalized harmonic function $\tilde \psi$ on $\Gamma =
\Gamma(A)$ defines a trace on $A$; in fact, a trace on each $A_n$
is determined by its value on minimal idempotents, so the
assignment
$$
\psi^{(n)}(e) = \tilde \psi(v),
\tag 3.7
$$
whenever $e$ is a minimal idempotent in $A_n$ such that $e M_v
\ne (0)$, defines a trace on $A_n$.  The harmonicity of $\tilde
\psi$ implies that the $\psi^{(n)}$ are {\it coherent}, i.e., the
restriction of $\psi^{(n+1)}$ from $A_{n+1}$ to the subalgebra
$A_n$ coincides with $\psi^{(n)}$. As a result, the traces
$\psi^{(n)}$ determine a trace of the limiting algebra
$A=\varinjlim A_n$.

The set of traces on $A$ is a compact convex set, with the
topology of pointwise convergence.  The map $\tilde \psi \mapsto
\psi$ is an affine homeomorphism between the space of positive
normalized harmonic functions on $\Gamma = \Gamma(A)$ and the
space of traces on $A$.  From the point of view of traces, the
Martin boundary of $\Gamma$ consists of traces $\psi$ which can
be obtained as limits of a sequence $\psi_n$, where $\psi_n$ is
an extreme trace on $A_n$.  All extreme traces on $A$ are in the
Martin boundary, so determination of the Martin boundary is a
step towards determining the set of extreme traces on $A$.

The locally semisimple algebra corresponding to the
Young-Fibonacci lattice $\yf$ is the Okada algebra $\Cal F$
introduced in Section 2.

\subhead The main result \endsubhead

We can now give a  description of the Martin boundary of the
Young-Fibonacci lattice  (and consequently of a  Poisson-type
integral representation for non-negative harmonic functions on
$\Bbb{YF}$).

\smallskip\noindent
{\bf Definition.}  Let $w$ be an infinite word in the alphabet
$\{1,2\}$ (infinite Fibonacci word), and let $d_1, d_2, \dots$
denote the positions of 2's in $w$.  The word $w$ is said to be
{\it summable} if, and only if, the series $\sum_{j=1}^\infty
1/d_j$ converges, or, equivalently, the product
$$
\pi(w) = \prod_{j:d_j\ge2} \Big(1 - {1 \over d_j}\Big) > 0
\tag 3.8
$$
converges.

\medskip
As for any differential poset, the lattice $\yf$ has  a
distinguished harmonic function $\varphi_P$, called the {\it
Plancherel harmonic function;} $\varphi_P$ is an element of the
Martin boundary. The complement of $\{\varphi_P\}$ in the Martin
boundary of $\Bbb{YF}$ can be  parametrized with two parameters
$(\beta, w)$; here $\beta$ is a real number, $0 < \beta \le 1$,
and $w$ is a summable infinite word in the alphabet $\{1,2\}$.

We denote by $\Omega$ the parameter space for the Martin
boundary:

\medskip
\noindent {\bf Definition.}
Let the space $\Omega$ be the union of a point $P$ and the set
$$
\{(\beta, w) : 0 < \beta \le 1, \   w \
\text{a summable infinite word in the alphabet  } \{1,2\}\},
$$
with the following topology:  A sequence $(\beta^{(n)}, w^{(n)})$
converges to $P$ iff
$$
\beta^{(n)} \rightarrow 0 \quad \text{or}\quad \pi(w^{(n)}) \rightarrow 0.
$$
A sequence $(\beta^{(n)}, w^{(n)})$ converges to $(\beta, w)$ if,
and only if,
$$
w^{(n)} \rightarrow w \ \text{\rm (digitwise)} \quad \text{and}
\quad \beta^{(n)} \pi(w^{(n)}) \rightarrow \beta \pi(w).
$$
\medskip

We will describe in Section 7 the mapping $\omega \mapsto
\varphi_\omega$ from $\Omega$ to the set of normalized positive
harmonic functions on $\yf$.

We are in a position now to state the main result of the paper.

\proclaim{Theorem 3.2}
The map $\omega\mapsto\varphi_\omega$ is a homeomorphism of the
space $\Omega$ onto the Martin boundary of the Young-Fibonacci
lattice. Consequently, for each  probability measure $M$  on
$\Omega$, the integral
$$
\varphi(v) = \int_\Omega \varphi_\omega(v)\, M(d\omega),\qquad v\in\Bbb{YF}
\tag 3.9
$$
provides a normalized, non-negative harmonic function on the
Young-Fibonacci lattice $\Bbb{YF}$. Conversely, every such
function admits an integral representation with respect to a
measure  $M$ on $\Omega$ (which may not be unique).
\endproclaim

In general, for all differential posets, we show that there is a
flow
$$
(t, \varphi) \mapsto C_t(\varphi)
$$
on $[0,1] \times \Cal H$ with the properties
$$
C_t(C_s(\varphi)) =
C_{ts}(\varphi) \quad \text{and} \quad C_0(\varphi) = \varphi_P.
\tag 3.10
$$
For the Young-Fibonacci lattice, one has $C_t(\varphi_{\beta, w})
= C_{t\beta, w}$ and $C_t(\varphi_P)=\varphi_P$. In particular,
the flow on $\Cal H$ preserves the Martin boundary. {\it It is
not clear whether this is a general phenomenon for differential
posets.}

We have not yet been able to characterize the extreme points
within the Martin boundary of $\yf$.  In a number of similar
examples, for instance the Young lattice, all elements of the
Martin boundary are extreme points.

\head \S 4. Harmonic functions on differential posets \endhead

The Young-Fibonacci lattice is an example of a differential
poset. In this section, we introduce some general constructions
for harmonic functions on a differential poset. Later on in
Section 7 we use the construction to obtain the Martin kernel of
the graph $\Bbb{YF}$.

\subhead Type I harmonic functions \endsubhead

In this subsection we don't need any special assumptions on the
branching diagram $\Gamma$. Consider an infinite path
$$
t = (v_0, v_1, \ldots, v_n, \ldots)
$$
in $\Gamma$.  For each vertex $u \in \Gamma$ the sequence $\{d(u,
v_n)\}_{n=1}^\infty$ is weakly increasing, and we shall use the
notation
$$
d(u,t) = \lim_{n \to \infty} d(u,v_n).
\tag 4.1
$$
Note that $d(u,t) = d(u,s)$ if the sequences $t,s$ coincide
eventually.

\proclaim{Lemma 4.1}
The following conditions are equivalent for a path $t$ in
$\Gamma$:
\smallskip
{\rm (i)}\quad All but finitely many vertices $v_n$ in
the path $t$ have a single immediate pre\-de\-ces\-sor $v_{n-1}$.
\smallskip
{\rm (ii)}\quad $d(\varnothing,t) < \infty$.
\smallskip
{\rm (iii)}\quad$d(u,t) < \infty$ for all $u \in \Gamma$.
\smallskip
{\rm (iv)}\quad There are only finitely many paths which
eventually coincide with $t$.
\endproclaim

\demo{Proof} It is clear that $d(\varnothing,v_{n-1}) =
d(\varnothing,v_n)$ iff $v_{n-1}$ is the only predecessor of
$v_n$. Since $d(u,t) < d(\varnothing,t)$ for all $u \in \Gamma$,
we have $(i) \to (ii) \to (iii) \to (i)$. The number of paths $s
\in T$, equivalent to $t$ is exactly $d(\varnothing,t)$.
\qed\enddemo

In case $\Gamma = \Y$ is the Young lattice, there are only
two paths (i.e. Young tableaux) satisfying these conditions: $t =
((1),\ldots,(n),\ldots)$ and $t = ((1),\ldots,(1^n),\ldots)$. In
case of Young-Fibonacci lattice there are countably many paths
satisfying the conditions of Lemma 4.1. The vertices of such a
path eventually take the form $v_n=1^{n-m}v$, $n\ge m$, for some
Fibonacci word $v$ of rank $m$. Hence, the equivalence class of
eventually coinciding paths in $\Bbb{YF}$ with the properties of
Lemma 4.1 can be labelled by infinite words in the alphabet
$\{2,1\}$ with only finite number of 2's. We denote the set of
such words as $1^\infty{\yf}$.

\proclaim{Proposition 4.2}
Assume  that a path $t $ in $\Gamma$ satisfies the conditions of
Lemma 4.1. Then
$$
\f_t(v) = {d(v,t) \over d(\varnothing,t)}, \quad v \in \Gamma
\tag 4.2
$$
is a positive normalized harmonic function on $\Gamma$.
\endproclaim

\demo{Proof} Since $d(v,t) = \sum_{w\colon v\nearrow w} d(w,t)$,
the function $\f_t$ is harmonic. Also, $\f_t(v) \ge 0$ for all $v
\in \Gamma$, and $\f_t(\varnothing) = 1$.
\qed\enddemo

We say that these harmonic functions  are {\it of type} I, since
the corresponding $AF$-algebra traces are  traces of finite -
dimensional irreducible representations (type I factor -
representations). It is clear that all the harmonic functions of
type I are indecomposable.

\subhead Plancherel harmonic function \endsubhead

Let us assume now that the poset $\Gamma$ is {\it differential}
in the sense of \cite{St1} or, equivalently, is a {\it self-dual
graph} in terms of \cite{F1}. The properties of differential
posets which we need are surveyed in the Appendix.

\proclaim{Proposition 4.3}
The function
$$
\f_P(v) = {d(\varnothing,v) \over n!}; \quad v \in \Gamma,\; n=|v|
\tag 4.3
$$
is a positive normalized harmonic function on the differential
poset  $\Gamma$.
\endproclaim

\demo{Proof}  This follows directly from (A.2.1) in the Appendix.
\qed\enddemo

Note that if $\Gamma = \Y$ is the Young lattice, the function
$\f_P$ corresponds to the Plancherel measure of the infinite
symmetric group (cf. \cite{KV}).

\subhead Contraction of harmonic functions on a differential
poset \endsubhead

Assume  that $\Gamma$ is a differential poset. We shall show that
for any harmonic function $\f$  there is a family of affine
transformations, with one real parameter $\tau$, connecting the
Plancherel function $\f_P$  to $\f$.

\proclaim{Proposition 4.4} For $0 \le \tau \le 1$ and a harmonic
function $\f$, define a function $C_\tau(\f)$ on the set of
vertices of the differential poset $\Gamma$ by the formula
$$
C_\t(\f)(v) = \sum_{k=0}^n {\tau^k(1-\tau)^{n-k} \over (n-k)!}
\sum_{|u|=k} \f(u) d(u,v), \quad n=|v|.
\tag 4.4
$$
Then $C_\t(\f)$ is a positive normalized harmonic function, and
the map $\varphi \mapsto C_t(\varphi)$ is affine.
\endproclaim

\demo{Proof} We introduce the notation
$$
S_k(v,\f) = \sum_{|u|=k} \f(u)\; d(u,v).
\tag 4.5
$$
First we observe the identity
$$
\sum_{w\colon v\nearrow w} S_k(w,\f) =
S_{k-1}(v,\f) + (n-k+1)\,S_k(v,\f),
$$
which is obtained from a straightforward computation using
(A.2.3) from the Appendix, and the harmonic property (3.1) of the
function $\f$. From this we derive that
$$
\aligned & \sum_{w\colon v\nearrow w} C_\t(\f)(w) =
\sum_{w\colon v\nearrow w} \sum_{k=0}^{n+1}
{\tau^k(1-\tau)^{n-k+1} \over (n-k+1)!} S_k(w,\f) = \\ = &
\sum_{k=1}^{n+1} {\tau^k(1-\tau)^{n-k+1} \over (n-k+1)!} S_{k-1}(v,\f) +
\sum_{k=0}^n {\tau^k(1-\tau)^{n-k+1} \over (n-k)!} S_k(v,\f) = \\
= & \tau \; C_\t(\f)(v) + (1-\tau)\; C_\t(\f)(v) = C_\t(\f)(v).
\endaligned
\tag 4.6
$$
This shows that $C_\t(\f)$ is harmonic. It is easy to see that
$C_\tau(\f)$ is normalized and positive, and that the map
$\varphi \mapsto C_t(\varphi)$ is affine.
\qed\enddemo

\medskip
\noindent {\bf Remarks.}
\noindent\quad {\rm (a)} The semigroup property holds:
$C_t(C_s(\varphi)) = C_{st}(\varphi)$;
\noindent\quad{\rm (b)} $C_0(\varphi) = \varphi_P$, for all
$\varphi$, and $C_t(\varphi_P) = \varphi_P$ for all $t$, $0\le t
\le 1$;
\noindent\quad {\rm (c)} $C_1(\varphi) = \varphi$.  These
statements can be verified by straightforward computations.

\medskip
\noindent {\bf Example.}
Let $\varphi$ denote the indecomposable harmonic function on the
Young lattice with the Thoma parameters $(\alpha;\beta;\gamma)$,
see \cite{KV} for definitions. Then the function
$C_\tau(\varphi)$ is also indecomposable, with the Thoma
parameters $(\tau\alpha;\tau\beta;1-\tau(1-\gamma))$.

\subhead Central measures and contractions\endsubhead

Recall (see \cite{KV}) that for any harmonic function $\f$ on
$\Gamma$ there is a {\it central measure} $M^\varphi$ on the
space $T$ of paths of $\Gamma$, determined by its level
distributions
$$
M^\varphi_n(v) = d(\varnothing,v)\, \f(v), \quad v \in \Gamma_n.
\tag 4.7
$$
In particular, $\sum_{v \in \Gamma_n} M^\varphi_n(v) = 1$ for all
$n$.

There is a simple probabilistic description of the central
measure corresponding to a harmonic function on a differential
poset obtained by the contraction of Proposition 4.4. Define a
random vertex $v \in \Gamma_n$ by the following procedure:

(a)\quad Choose a random $k$, $0 \le k \le n$ with the binomial
distribution
$$
b(k) = {n \choose k} \tau^k (1-\tau)^{n-k}
\tag 4.8
$$

(b)\quad Choose a random vertex $u \in \Gamma_k$ with probability
$$
M_k^{\f}(u) = {d(\varnothing,u) \f(u)}
\tag 4.9
$$

(c)\quad Start a random walk at the vertex $u$, with the Plancherel
transition probabilities
$$
p_{x,y} = {d(\varnothing,y) \over (r+1)\, d(\varnothing,x)}; \quad
|x|=r, \; x \nearrow y.
\tag 4.10
$$
Let $v$ denote the vertex at which the random walk first hits the
$n$'th level set $\Gamma_n$. We denote by $M_n^{(\tau,\f)}$ the
distribution of the random vertex $v$.

\proclaim{Proposition 4.5}
The distribution $M_n^{(\tau,\f)}$ is the $n$'th level
distribution of the central measure corresponding to the
harmonic function $C_\t(\f)$:
$$
M_n^{(\tau,\f)}(v) = d(\varnothing,v)\, C_\t(\f)(v).
\tag 4.11
$$
\endproclaim

\demo{Proof} It follows from (A.2.1) that (4.10) is a probability
distribution. By Lemma A.3.2, the probability to hit $\Gamma_n$
at the vertex $v$, starting the Plancherel walk at $u \in
\Gamma_k$, is
$$
p(u,v) = {k! \over n!}\, {d(u,v)\, d(\varnothing,v) \over
d(\varnothing,u)}.
\tag 4.12
$$
The Proposition now follows from the definition of the
contraction $C_\t(\f)$ written in the form
$$
C_\t(\f)(v)\, d(\varnothing,v) =
\sum_{k=0}^n b(k) \sum_{u \in \Gamma_k} M_k^{\f}(u)\, p(u,v).\qed
\tag 4.13
$$
\enddemo

\noindent {\bf Example.}
Let $\Gamma = \Y$ be the Young lattice and let $t =
((1),(2),\ldots,(n),\ldots)$ be the one-row Young tableau. Then
the distribution (4.9) is trivial, and the procedure reduces to
choosing a random row diagram $(k)$ with the distribution (4.8)
and applying the Plancherel growth process until the diagram
gains $n$ boxes.

\head \S 5. Okada clone of the symmetric function ring  \endhead

In this Section we introduce the Okada variant of the symmetric
function algebra, and its two bases analogous to the Schur
function basis and the power sum basis. The Young-Fibonacci
lattice arises in a Pieri-type formula for the first basis.

\subhead  The rings $R$ and $R_\infty$ \endsubhead

Let $R = \Bbb Q<X,Y>$ denote the ring of all polynomials in two
non-commuting variables $X,Y$. We endow $R$ with a structure of
graded ring, $R = \bigoplus_{n=0}^\infty R_n$, by declaring the
degrees of variables to be $\deg X = 1$, $\deg Y = 2$.  For each
word
$$
v = 1^{k_t}\, 21^{k_{t-1}}\, \ldots 1^{k_1}21^{k_0} \in \yf_n,
\tag 5.1
$$
let $h_v$ denote the monomial
$$
h_v = X^{k_0}YX^{k_1}\, \ldots X^{k_{t-1}}Y\, X^{k_t}
\tag 5.2
$$
Then $R_n$ is a $\Q$-vector space with the $f_n$ (Fibonacci
number) monomials $h_v$ as a basis.

We let $R_\infty=\varinjlim R_n$ denote the inductive limit of
linear spaces $R_n$, with respect to imbeddings $Q\mapsto Q\,X$.
Equivalently, $R_\infty=R/(X-1)$ is the quotient of $R$ by the
principal left ideal generated by $X-1$.  Linear functionals on
$R_\infty$ are identified with linear functionals $\varphi$ on
$R$ which satisfy $\varphi(f) = \varphi(f X)$. The ring
$R_\infty$ has a similar r\^ole for the Young-Fibonacci lattice
and the Okada algebra $\Cal F$ as the ring of symmetric functions
has for the Young lattice and the group algebra of the infinite
symmetric group $\frak{S}_\infty$ (see [M]).

\subhead  Non-commutative Jacobi determinants  \endsubhead

The following definition is based on a remark which appeared in
the preprint version of \cite{Ok}. We consider two
non-commutative $n$-th order determinants
$$
P_n = \left| \matrix
X & Y & 0 & 0 & \ldots & 0 & 0\\
1 & X & Y & 0 & \ldots & 0 & 0\\
0 & 1 & X & Y & \ldots & 0 & 0\\
  &   &   & \ldots &   &   &  \\
0 & 0 & 0 & 0 & \ldots & 1 & X
\endmatrix \right|
\tag 5.3
$$
and
$$
Q_{n-1} = \left| \matrix
Y & Y & 0 & 0 & \ldots & 0 & 0\\
X & X & Y & 0 & \ldots & 0 & 0\\
0 & 1 & X & Y & \ldots & 0 & 0\\
  &   &   & \ldots &   &   &  \\
0 & 0 & 0 & 0 & \ldots & 1 & X
\endmatrix \right|.
\tag 5.4
$$
By definition, the non-commutative determinant is the expression
$$
\det (a_{ij}) = \sum_{w \in \frak S_n} \text{\rm sign}(w)\,
a_{w(1)1}\,a_{w(2)2}\, \ldots a_{w(n)n}.
\tag 5.5
$$
In other words, the $k$-th factor of every summand is taken from
the $k$-th column. Note that polynomials (5.3), (5.4) are
homogeneous elements of $R$, $\deg P_n=n$ and $\deg Q_{n-1}=n+1$.

Following Okada, we define {\it Okada-Schur polynomials} (or
$s$-functions) as the products
$$
s_v = P_{k_0}\, Q_{k_1}\, \ldots Q_{k_t}, \qquad
v=\underbrace{1^{k_t}2}\,\ldots\underbrace{1^{k_1}2}\,1^{k_0}\in\yf_n
\tag 5.6
$$
(cf. \cite{Ok}, Proposition 3.5). The polynomials $s_v$ for $|v|
= n$  are homogeneous of degree $n$, and form a basis of  the
linear space $R_n$.
We define a scalar product $\langle\,.\,,.\,\rangle$ on the space $R$
by declaring $s$-basis to be orthonormal.

\subhead The branching of Okada-Schur functions \endsubhead

We will use the formulae
$$
P_{n+1} = P_n\,X - P_{n-1}\,Y
\quad n\ge 1,
\tag 5.7
$$
$$
Q_{n+1} = Q_n\,X - Q_{n-1}\,Y,
\quad n\ge 1,
\tag 5.8
$$
obtained by decomposing the determinants (5.3), (5.4) along the
last column. The first identity is also true for $n=0$, assuming
$P_{-1}=0$. The $n=0$ case of the second identity (5.8) can be
written in the form
$$
Q_0\, X = X\, Q_0 + Q_1.
\tag 5.9
$$
One can think of (5.9) as of a commutation rule for passing $X$ over a
factor of type $Q_0$. It is clear from (5.9) that
$$
Q_0^m\, X = X\, Q_0^m + \sum_{i=1}^m Q_0^{m-i}\,Q_1\,Q_0^{i-1}.
\tag 5.10
$$
It will be convenient to rewrite (5.7), (5.8) in a
form similar to (5.9):
$$
P_n\, X = P_{n+1} + P_{n-1}\,Q_0.
\tag 5.11
$$
$$
Q_n\, X = Q_{n+1} + Q_{n-1}\,Q_0,
\tag 5.12
$$
The following formulae are  direct consequences of (5.10) -- (5.12):
$$
P_n\, Q_0^m \, X = \sum_{i=0}^m P_n\, Q_0^{m-i}\, Q_1\, Q_0^i +
P_{n+1}\, Q_0^m + P_{n-1}\, Q_0^{m+1},
\tag 5.13
$$
$$
Q_n\, Q_0^m \, X = \sum_{i=0}^m Q_n\, Q_0^{m-i}\, Q_1\, Q_0^i +
Q_{n+1}\, Q_0^m + Q_{n-1}\, Q_0^{m+1}.
\tag 5.14
$$
It is understood in (5.13), (5.14) that $n\ge 1$.

The formulae (5.10), (5.13) and (5.14) imply

\proclaim{Theorem 5.1 \rm(Okada)}  For every  $w \in \yf_n$ the
product of the Okada-Schur determinant $s_w$ by $X$ from the
right hand side can be written as
$$
s_w\, X = \sum_{v: w \nearrow v} s_v.
\tag 5.15
$$
\endproclaim

This theorem says that the branching of Okada $s$-functions
reproduces the branching law for the Young-Fibonacci lattice.  In
the following statement, $U$ is the ``creation operator" on
$\Fun(\yf)$, which is defined in the Appendix, (A.1.1).

\proclaim{Corollary 5.2} The assignment $\Theta : v \mapsto s_v$
extends to a linear isomorphism
$$
\Theta : \cup_n \Fun(\yf_n) \rightarrow R
$$
taking $\Fun(\yf_n)$ to $R_n$ and satisfying $\Theta \circ U(f) =
\Theta(f) X$.
\endproclaim

Because of this, we will sometimes write $U(f)$ instead of $f X$
for $f \in R$, and $D(f)$ for $\Theta\circ D\circ
\Theta^{-1}(f)$, see (A.1.2) for the definition of $D$.

\proclaim{Corollary 5.3} There exist one-to-one correspondences between:
\smallskip\noindent
{\rm(a)}\quad Non-negative, normalized harmonic functions on $\yf$;
\smallskip\noindent
{\rm(b)}\quad Linear functionals $\varphi$ on $\Fun(\yf)$ satisfying
$$
\varphi\circ U = \varphi,\quad  \varphi(1) = 1,\quad \text{ and} \quad
\varphi(\delta_v) \ge 0 \qquad \text{ for } v \in \yf;
$$
\smallskip\noindent
{\rm(c)}\quad Linear functionals $\varphi$ on $R$ satisfying
$$
\varphi(f) = \varphi(fX) \quad \text{ for all } f \in R,
\quad \varphi(1) = 1\quad\text{and}\quad \varphi(s_v) \ge 0,
\qquad \text{ for } v \in \yf;
$$
\smallskip\noindent
{\rm(d)}\quad Linear functionals on $R_\infty = \varinjlim R_n$ satisfying
$$
\varphi(1) = 1\quad\text{and}\quad \varphi(\bar s_v) \ge 0,
\qquad \text{ for } v \in \yf,
$$
where $\bar s_v$ denotes the image of $s_v$ in $R_\infty$;
\smallskip\noindent
{\rm(e)}\quad  Traces of the Okada algebra $\Cal F_\infty$.
\endproclaim

We refer to linear functionals $\varphi$ on $R$ satisfying
$\varphi(s_v) \ge 0$ as {\it positive} linear func\-tio\-nals.

\subhead The Okada $p$-functions \endsubhead

Following Okada \cite{Ok}, we introduce another family of
homogeneous polynomials, labelled by Fibonacci words $v\in\yf$,
$$
p_v = \big(X^{k_0+2}-(k_0+2)X^{k_0}Y\big)\ldots
\big(X^{k_{t-1}+2}-(k_{t-1}+2)X^{k_{t-1}}Y\big)\,X^{k_t},
\tag 5.16
$$
where
$$
v=1^{k_t}\underbrace{21^{k_{t-1}}}\,\ldots\underbrace{21^{k_0}}.
$$
One can check that $\{p_v\}_{|v|=n}$ is a $\Bbb Q$-basis of
$R_n$.  Two important properties of the $p$-basis which were
found by Okada are:
$$
U(p_v) = p_v X = p_{1v} \quad \text{and} \quad D(p_{2v}) = 0. \tag 5.17
$$
Since the images of $p_u$ and of $p_{1u}$ in $R_\infty$ are the
same, we can conveniently denote the image  by $p_{1^\infty u}$.
The family of $p_v$, where $v$ ranges over $1^\infty \yf$, is a
basis of $R_\infty$.

\subhead  Transition matrix from $s$-basis to $p$-basis \endsubhead

We  denote the transition matrix relating the two bases $\{p_u\}$
and $\{s_v\}$ by $X^v_u$; thus
$$
p_u = \sum_{|v|=n} X_u^v\, s_v, \quad u, v  \in  \yf_n.
\tag 5.18
$$
The coefficients $X_u^v$ are analogous to the character matrix of
the symmetric group $\frak S_n$. They were described recurrently
in \cite{Ok}, Section 5, as follows:
$$
X_{1u}^{1v}=X_u^v;\quad X_{2u}^{1v}=X_{1u}^v;\quad X_{2u}^{2v}=-X_u^v,
\tag 5.19
$$
$$
X_{11u}^{2v} = (m(u)+1)\, X_u^v;\quad X_{12u}^{2v} = 0,
\tag 5.20
$$
where $m(u)$ is defined in (6.2) below. An  explicit product
expression for the $X^v_u$ will be given in the next section.

\head \S 6. A product formula for Okada characters \endhead

In this Section we improve Okada's description of the character
matrix $X_u^v$ to obtain the product formula (6.11) and its
consequences.

\subhead Some notation \endsubhead

We recall some notation from \cite{Ok} which will be used below.
Let $v$ be a  Fibonacci word:
$$
v = 1^{k_t}\, 21^{k_{t-1}}\,\ldots 1^{k_1}21^{k_0} \in \yf_n.
$$
Then:

\noindent (6.1)\quad $\epsilon(v)=+1$ if the  rightmost digit of
$v$ is $1$, and $\epsilon(v)=-1$ otherwise.

\noindent (6.2)\quad $m(v)=k_t$ is the number of $1$'s at the
left end of $v$.

\noindent(6.3)\quad The rank of $v$, denoted $|v|$, is the sum of
the digits of $v$.

\noindent(6.4) \quad  If $v = v_1 2 v_2$, then the {\it position}
of the indicated 2 is $|v_2| + 1$.

\noindent(6.5)\quad $d(v)=\prod_{i=0}^{t-1}
(k_0+\ldots+k_i+2i+1)$. In other words, $d(v)$ is the product
of the positions of 2's in $v$.

\noindent(6.6)\quad $z(v) = k_t!\, (k_{t-1}+2)k_{t-1}!\, \ldots
(k_0+2)k_0!$.

\noindent(6.7)\quad The {\it block ranks} of $v$ are the numbers
$k_0+2, k_1+2, \dots, k_{t-1}+2, k_t$.

\noindent(6.8)\quad The {\it inverse block ranks} of $v$ are
$k_t+2, k_{t-1}+2, \dots, k_1+2, k_0$.

Consider a sequence $\bar{n}=(n_t,\ldots,n_1,n_0)$ of positive
integers with $\sum n_i = n$. We call a word $v \in \yf_n$ {\it
$\bar n$-splittable}, if it can be written as a concatenation
$$
v = v_t\, \ldots v_1\,v_0, \quad \text{where } |v_i|=n_i \quad
\text{\rm for } i=0,1,\ldots,t.
\tag 6.8
$$

\proclaim{Lemma 6.1}
Let $\bar{n}=(n_t,\ldots,n_1,n_0)$ be the sequence of block ranks
in a Fibonacci word $u$. Then

(i)\quad $X_u^v \ne 0$ if, and only if, the word $v$ is $\bar
n$-splittable.

(ii)\quad If $v=v_t\ldots v_1 v_t$ is the $\bar n$-splitting,
then
$$
X_u^v = d(v_t) g(v_{t-1}) \ldots g(v_0),
\tag 6.9
$$
where
$$
g(w) = \cases
+d(w'), & \text{ if } w=1w'\\
-d(w'), & \text{ if } w=2w'.
\endcases
\tag 6.10
$$
\endproclaim
\demo{Proof} This is a direct consequence of Okada's recurrence
relations cited in the previous section.\qed \enddemo

\proclaim{Proposition 6.2}
Let $u, v \in \yf_n$. Let $\delta_1, \delta_2, \dots, \delta_m$
be the positions of 2's in the word $u$,  and put $\delta_{m+1} =
\infty$.  Let $d_1, d_2, \dots, d_r$ the positions of 2's in
the word $v$.  Then
$$
X_u^v = \prod_{j = 1}^m \prod_{\delta_j \le d_s <
\delta_{j+1}}\big(d_s - (\delta_j + 1) \big).
\tag 6.11
$$
\endproclaim

\demo{Proof}  This can also be derived directly from Okada's
recurrence relations, or from the previous lemma.  Note in
particular that $X_u^v = 0$ if, and only if, $d_s = \delta_j+1$
for some $s$ and $j$; this is the case if, and only if, $v$ does
not split according to the block ranks of $u$.
\qed \enddemo

We define $\tilde X_u^v = {d(v)}^{-1}X^v_u$; from Proposition 6.2
and the dimension formula (6.5), we have the expression
$$
\tilde X_u^v = \prod_{j = 1}^m \prod_{\delta_j \le d_s <
\delta_{j+1}}\left(1 - \frac{\delta_j + 1}{d_s} \right)
\tag 6.12
$$

\subhead  The inverse transition matrix  \endsubhead

According to [Ok], Proposition 5.3, the inverse formula to
Equation (5.18) can be written in the form
$$
s_v = \sum_{|u|=n} X_u^v\, {p_u \over z(u)}, \quad v \in \yf_n,
\tag 6.13
$$
We will give  a description of a  column $X_u^v$ for a fixed $v$.

\proclaim{Lemma 6.3}
Let $\bar n = (n_t, \ldots, n_1,n_0)$ be the sequence of  inverse
block ranks $n_t=k_t+2, \ldots, n_1=k_1+2, n_0=k_0$ in a word $v
= (1^{k_t}2\, \ldots 1^{k_1}2\, 1^{k_0}) \in \yf_n$. Then

(i)\quad $X_u^v \ne 0$ if, and only if,  the word $u$ is $\bar
n$-splittable.

(ii)\quad If $u=u_t \ldots u_1 u_0$ is the $\bar n$-splitting for
$u$, then
$$
X_u^v = f_1\, f_2\, \ldots f_t,
\tag 6.14
$$
where
$$
f_j = \cases -1, & \text{ if } \epsilon(u_j)=-1\\
1+m(u_{j-1} \ldots u_1 u_0), & \text{ if } \epsilon(u_j)=+1.
\endcases
\tag 6.15
$$
Here $m(u)=m$ denotes the number of $1$'s at the left end of
$u=1^m2 u'$.
\endproclaim
\demo{Proof} This is another corollary of Okada's recurrence
relations cited in Section 5.
\qed\enddemo

\head \S 7.  The Martin Boundary of the Young-Fibonacci Lattice
\endhead

In this section, we examine certain elements of the Martin
boundary of the Young-Fibonacci lattice $\Bbb{YF}$.  Ultimately
we will show that the harmonic functions listed here comprise the
entire Martin boundary.

It will be  useful for us to evaluate normalized positive linear
functionals on the ring $R_\infty$  (corresponding to normalized
positive harmonic functions on $\yf$) on the basis $\{p_u\}$.
The first result in this direction is the evaluation of the
Plancherel functional on these basis elements.

\proclaim{Proposition 7.1} $\varphi_P(p_u) = 0$ for all
Fibonacci words $u$ containing at least one 2.
\endproclaim

\demo{Proof} It follows from the definition of the Plancherel
harmonic function $\varphi_P$ that for $w \in \yf_n$,
$$
\sum_{v: v \nearrow w} \varphi_P(v) = n\; \varphi_P(w).
$$
Therefore, for all $f \in R_n$,
$$
\varphi_P(D f) = n\; \varphi_P(f).
$$
If $u = 1^\infty 2 v$, and $|2v|=n$, then
$$
\varphi_P(p_u) = \varphi_P(p_{2v}) =
\frac{1}{n}\; \varphi_P(D p_{2v}) = 0,
$$
since $Dp_{2v} = 0$, by [Ok], Proposition 5.4.
\qed\enddemo

For each word $w \in \yf_n$,  the path $(w, 1w, 1^2w, \dots)$
clearly satisfies the conditions of Proposition 4.1, and
therefore there is a type I harmonic function on $\yf$ defined by
$$\psi_w(v) = \lim_{k\rightarrow \infty} \frac{d(v, 1^k
w)}{d(\varnothing, 1^k w)}.$$

\proclaim{Proposition 7.2}
Let $w \in \yf_n$, and let $d_1, d_2, \dots$ be the positions of
2's in $w$. Let $u$ be a word in $1^\infty \yf$ containing at
least one 2, and let $\delta_1, \delta_2, \dots, \delta_m$ be the
positions of 2's in $u$.  Then:
$$
\psi_w(p_u) = \prod_{i=1}^m \prod_{\delta_i\le d_j<\delta_{i+1}}
\left(1 - {\delta_i + 1 \over d_j}\right).
\tag 7.1
$$
\endproclaim

\demo{Proof} Let $u = 1^\infty u_0$, where $u_0 \in \yf_m$.
Choose $r, s \ge 1$ such that $|1^s w| = |1^r u_0|$.  Then
$$
\aligned
\psi_w(p_u)
&=(d(1^s w))^{-1}\langle p_{1^r u_0}, s_{1^s w} \rangle \cr&=
(d(1^s w))^{-1}\langle \sum_v X^v_{1^r u_0} s_v, s_{1^s w} \rangle \cr &=
(d(1^s w))^{-1}X^{1^s w}_{1^ru_0}.
\endaligned
$$
Thus the result follows from Equation (6.12).
\qed\enddemo

Next we describe some harmonic functions which arise from
summable infinite  words. Given a summable infinite word $w$,
define a linear functional on the ring $R_\infty$  by the
requirements $\varphi_w(1)=1$ and
$$
\varphi_w(p_u) =
\prod_{i=1}^m \prod_{\delta_i\le d_j<\delta_{i+1}}
\left(1 - {\delta_i + 1 \over d_j}\right),
\tag 7.3
$$
where $u \in 1^\infty \yf$. As usual, $\delta_1, \dots, \delta_m$
are the positions of 2's in $u$, and the $d_j$'s are the
positions of 2's in $w$.  It is evident that $\varphi_w(p_uX) =
\varphi_w(p_{1u}) = \varphi_w(p_u)$, so that $\varphi_w$  is in
fact a functional on $R_\infty$.

\proclaim{Proposition 7.3}
If $w$ is a summable infinite Fibonacci word, then $\varphi_w$ is
a normalized positive linear functional on $R_\infty$, so
corresponds to a normalized positive harmonic function on $\yf$.
\endproclaim

\demo{Proof} Only the positivity needs to be verified. Let $w_n$
be the finite word consisting of the rightmost $n$ digits of $w$.
It follows from the product formula for the normalized characters
$\psi_{w_n}$ that $\varphi_w(p_u) =
\displaystyle\lim_{n\rightarrow \infty} \psi_{w_n}(p_u) $.
Therefore also $\varphi_w(s_v) =
\displaystyle\lim_{n\rightarrow\infty} \psi_{w_n}(s_v) \ge 0 $.
\qed\enddemo

Given a summable infinite Fibonacci word $w$ and $0 \le \beta \le
1$, we can define the harmonic function $\varphi_{\beta, w}$ by
contraction of $\varphi_w$, namely, $\varphi_{\beta, w} =
C_\beta(\varphi_w)$.

For $u \in 1^\infty \yf$, we let $||u||$ denote the {\it
essential rank} of $u$, namely $||u|| = 1 + \delta$, where
$\delta$ is the position of the leftmost 2 in $u$, and $||u|| =
0$ for $u = 1^\infty$.

\proclaim{Proposition 7.4}
Let $w$ be a summable infinite word and $0 \le \beta \le 1$. Let
$u \in 1^\infty \yf$.  Then
$$
\f_{\beta,w}(p_u) = \beta^{||u||}\f_w(p_u).
\tag 7.4
$$
\endproclaim

\demo{Proof}
The case $u = 1^\infty$ is trivial. Let $u = 1^\infty u_0$, where
$||u|| = |u_0| = n > 0$. Then for any linear functional
$\varphi$ on $R_\infty$,  one has $\varphi(p_u) =
\varphi(p_{u_0})$. For $0\le k \le n$, define
$$
A^w_{k,n} = \sum_{|v| = k} \sum_{|x| = n}\varphi_w(v) d(v,x) s_x.
$$
In particular,
$$
A^w_{n,n} = \sum_{|v| = n} \varphi_w(v) s_v.
$$
Note that $U A^w_{k,n-1} = A^w_{k,n}$, when $k\le n-1$.  It
follows from the definitions of $\f_{\beta,w}$ (cf. (4.4)) and of
$A^w_{k,n}$ that
$$
\f_{\beta,w}(f) =
\langle \sum_{k=0}^n {\beta^k(1-\beta)^{n-k} \over (n-k)!}
A^w_{k,n}, f \rangle,
$$
for $f$ in $R_n$, where $\langle \cdot, \cdot \rangle$ denotes
the inner product on $R$ with respect to which the Okada
$s$-functions form an orthonormal basis.  Recall that the
operators $U$ and $D$ are conjugate with respect to this inner
product.  Consequently,
$$
\aligned
\f_{\beta,w}(p_{u_0}) &= \langle \sum_{k=0}^n {\beta^k(1-\beta)^{n-k} \over
(n-k)!} A^w_{k,n}, p_{u_0} \rangle \cr
&= \langle \sum_{k=0}^{n-1} {\beta^k(1-\beta)^{n-k} \over (n-k)!}
U A^w_{k,n-1}, p_{u_0} \rangle
+ \langle\beta^n A^w_{n,n}, p_{u_0} \rangle \cr
&= \langle\beta^n A^w_{n,n}, p_{u_0} \rangle = \beta^n \f_w(p_{u_0}),
\endaligned
$$
since $\langle U A^w_{k,n-1}, p_{u_0} \rangle = \langle
A^w_{k,n-1}, D p_{u_0} \rangle = 0$ and $\langle  A^w_{n,n}, f
\rangle = \f_w( f)$ for $f$ in  $R_n$. \qed
\enddemo

\proclaim{Corollary 7.5}
The functionals $\varphi_{\beta, w}$ for $\beta > 0$ and $w$
summable are pairwise distinct, and different from the Plancherel
functional $\varphi_P$.
\endproclaim

\demo{Proof}
Suppose that $w$ is a summable word and that $A$ is the set of
positions of 2's in $w$. We set
$$
\pi_k(w) = \prod_{j: d_j \ge k-2}  \left(1- \frac{k}{d_j}\right).
$$
Then for each $k \ge 2$  and $\beta>0$,
$$
\varphi_{\beta, w}(p_{21^{k-1}}) = \beta^k \pi_k(w)
$$
is zero if and only if $k \in A$.  In particular $\varphi_{\beta,
w} \ne \varphi_P$, by Lemma 7.1, and moreover, $A\setminus\{1\}$
is determined by the sequence of values
$\varphi_{\beta,w}(p_{21^{k-2}})$, $k\ge2$. It is also clear that
$\varphi_{\beta,w}(p_2)=\beta^2\,\pi_2(w)$ is negative iff $1\in
A$, hence the set $A$ and therefore also $w$ are determined by
the values $\varphi_{\beta,w}(p_{21^k})$, $k\ge0$. Finally,
$\beta$ is determined by
$$
\beta =
\left(\frac{\varphi_{\beta,w}(p_{21^{k-1}})}{\pi_k(w)}\right)^{1/k},
$$
for any $k \not\in A$.
\qed\enddemo

\proclaim{Proposition 7.6}
For each summable infinite Fibonacci word $w$ and each $\beta$,
$0 \le \beta \le 1$, there exists a sequence $v^{(n)}$ of
finite Fibonacci words such that $\varphi_{\beta, w} =
\displaystyle\lim_{n\rightarrow \infty} \psi_{v^{(n)}}$.
\endproclaim

\demo{Proof}
If $\beta = 1$, put $r_n =0$; if $\beta = 0$, put $r_n = n^2$;
and if $0 < \beta < 1$, choose the sequence $r_n$ so that
$$
\lim_{n\rightarrow \infty} \frac{r_n}{n} = \frac{1-\beta^2}{\beta^2}.
$$
Then, in every case,
$$
\lim_{n\rightarrow \infty} \frac{n}{n+ r_n} = \beta^2.
$$
Let $w_n$ be the finite word consisting of the rightmost $n$ digits of
$w$, put $s_n = 2n+1-|w_n|$,  and
$$
v^{(n)} =  2^{r_n}1^{s_n}w_n.
$$
Fix $u = 1^\infty u_0 \in 1^\infty \yf$ and let $n \ge |u_0|$.
Suppose that $u_0$ has 2's at positions
$\delta_1,\delta_2,\dots,\delta_m$, and put $k = ||u|| = \delta_m
+ 1$.  Using the product formula for $\psi_{v^{(n)}}$, one
obtains
$$
\psi_{v^{(n)}}(p_u) = \psi_{w_n}(p_u) \left[(1-\frac{k}{2n+2})(1-
\frac{k}{2n+4})\cdots
(1-\frac{k}{2n+2r_n})\right].
$$
The first factor converges to $\varphi_w(p_u)$, so it suffices,
by Proposition 7.4, to show that the second factor converges to
$\beta^k$. The second factor reduces to
$$
\frac{\Gamma(n+r_n-k/2 + 1)
\Gamma(n+2)}{\Gamma(n+2-k/2)\Gamma(n+r_n + 1)}.
$$
Using the well-known fact that
$$
\lim_{n\rightarrow \infty} n^{b-a}\frac{\Gamma(n+a)}{\Gamma(n+b)} = 1,
$$
one obtains that the ratio of gamma functions is asymptotic to
$$
\left(\frac{n+2}{n+r_n} \right)^{k/2},
$$
which, by our choice of $r_n$, converges to $\beta^k$, as desired.
\qed\enddemo

This proposition shows that $\varphi_P$ as well as all of the
harmonic functions $\varphi_{\beta, w}$ are contained in the
Martin boundary of the Young-Fibonacci lattice.  In the following
sections, we will show that these harmonic functions make up the
entire  Martin boundary.

\head \S 8. Regularity conditions \endhead

In this Section we obtain a simple criterion for a sequence of
characters of finite dimensional Okada algebras to converge to a
character of the limiting infinite dimensional algebra $\Cal
F=\varinjlim\Cal F_n$. Using this criterion, the {\it regularity
conditions}, we show that the harmonic functions provided by the
formulae (7.3) and (7.4) make up the entire Martin boundary of
the Young-Fibonacci graph $\Bbb{YF}$. Technically, it is more
convenient to work with linear functionals on the spaces $R_n$
and their limits in $R_\infty=\varinjlim R_n$, rather than with
traces on $\Cal F$. As it was already explained in Section 5,
there is a natural one-to-one correspondence between traces of
Okada algebra $\Cal F_n$, and positive linear functionals on the
space $R_n$.

In this Section we shall use the following elementary inequalities:
$$
\Big(1 - {k\over d}\Big) \le
\Big(1 - {1\over d}\Big)^k,
\tag 8.1
$$
for every pair of positive integer numbers $d\ge2$ and $k$;
$$
\Big(1 - {1\over d}\Big)^{k^2} \le
\Big(1 - {k\over d}\Big),
\tag 8.2
$$
and
$$
1 \le \frac{\Big(1 - {1\over d}\Big)^k}{\Big(1 - {k\over d}\Big)} \le
1 + \frac{{k \choose 2}}{(d-k)^2},
\tag 8.3
$$
for every pair of integers $d>k$. We omit the straightforward
proofs of these inequalities.

\subhead  Convergence to the Plancherel measure \endsubhead

We first examine the important particular case of convergence to
the Plancherel character $\varphi_P$.

\smallskip\noindent
{\bf Definition.}
We define the function $\pi$ of a finite or summable word $v$ by
$$
\pi(v) = \prod_{j:\,d_j\ge2}\, (1 - \frac{1}{d_j}),
$$
where the $d_j$ are the positions of the 2's in $v$. We also
recall that for each $k\ge2$ the function $\pi_k$ was defined as
$$
\pi_k(v) = \prod_{j\, :\, d_j \ge k-1}\, (1 - \frac{k}{d_j}).
$$
Note that if $u=1^\infty21^{k-2}$, and $v$ is a summable word,
then $\varphi_v(p_u) = \pi_k(v)$, according to Equation (7.3).

\proclaim{Proposition 8.1}
The following properties of a sequence $w_n\in \yf$,
$n=1,2,\ldots$, are equivalent:

\noindent
{\rm(i)} The normalized characters $\psi_{w_n}$ converge to the
Plancherel character, i.e.,
$$
\lim_{n\to\infty} \psi_{w_n}(p_u) = \varphi_P(p_u),
$$
for each $u\in 1^\infty \yf$;
$$
\lim_{n\to\infty} \pi_k(w_n) = 0;\,
\tag ii
$$
for every $k=2,3,\ldots$;
$$
\lim_{n\to\infty} \pi(w_n) = 0.
\tag iii
$$
\endproclaim

The proof is based on the following lemmas.

\proclaim{Lemma 8.2}
For every finite word $v\in \yf$, and for every $u\in1^\infty\yf$
of essential rank $k=||u||$,
$$
|\psi_v(p_u)| \le |\pi_k(v)|.
\tag 8.4
$$
\endproclaim

\demo{Proof}
Let $\delta_1,\ldots,\delta_m$ indicate the positions of $2$'s in
the word $u$, and let $d_1,\ldots,d_n$ be the positions of $2$'s
in $v$. The essential rank of $u$ can be written as
$k=||u||=\delta_m+1$, so that
$$
\pi_k(v) = \prod_{j:d_j\ge\delta_m} \big(1 - {k\over d_j}\big).
$$
By the product formula,
$$
|\psi_v(p_u)| = |\pi_k(v)| \prod_{i=1}^{m-1} \prod_{\delta_i\le
d_j<\delta_{i+1}} \big|1 - {\delta_i + 1 \over d_j}\big| \le
|\pi_k(v)|,
$$
since none of  the factors in the product  exceed $1$. In fact,
$|1-(\delta+1)/d|=1/\delta\le1$ if $d=\delta$, $1-(\delta+1)/d=0$
if $d=\delta+1$, and $0\le1-(\delta+1)/d<1$ if $d>\delta+1$.
\qed\enddemo

\proclaim{Lemma 8.3}
For each $k=2,3,\ldots$, and for every word $v\in \yf$,
$$
|\pi_k(v)| \le \big(k\; \pi(v)\big)^k.
$$
\endproclaim

\demo{Proof}
It follows from (8.1) that
$$
\gathered
|\pi_k(v)| = \prod_{d_j\ge k-1} \big|1 - {k\over d_j}\big| \le
\prod_{d_j\ge k+1} \big(1 - {k\over d_j}\big) \le \\ \le
\prod_{d_j\ge k+1} \big(1 - {1\over d_j}\big)^k =
\pi^k(v)\; \prod_{2\le d_j\le k} \big(1 - {1\over d_j}\big)^{-k}.
\endgathered
$$
Since $(1-1/d)^{-1}>1$ for $d\ge2$, the last product can be
estimated as
$$
\prod_{2\le d_j\le k} \big(1 - {1\over d_j}\big)^{-k} \le
\big({1\over2}\;{2\over3}\; \ldots\; {k-1\over k}\big)^{-k} = k^k,
$$
and the lemma follows.
\qed\enddemo

\proclaim{Lemma 8.4}
If $d_1(v)\ne2$, then $|\pi_2(v)|\ge\pi^4(v)$, and if $d_1(v)=2$, then
$|\pi_3(v)|\ge\pi^9(v)$.
\endproclaim

\demo{Proof}
We apply the inequality (8.2). If $d_1(v)\ge3$, then
$$
|\pi_2(v)| = \prod_{j:d_j\ge3} \big(1 - {2 \over d_j}\big) \ge
\prod_{j:d_j\ge3} \big(1 - {1 \over d_j}\big)^4 = \pi^4(v)
$$
by the inequality (8.2). If $d_1(v)=1$, then $(1-2/d_1)=-1$, and
since $d_2\ge3$, the inequality holds in this case as well.

In case of $d_1(v)=2$ we have $d_2(v)\ge4$, hence
$$
|\pi_3(v)| = {1\over2} \prod_{j:d_j\ge4} \big(1 - {3\over d_j}\big); \qquad
|\pi(v)| = {1\over2} \prod_{j:d_j\ge4} \big(1 - {1\over d_j}\big),
$$
so that the second inequality of Lemma also follows from (8.2).
\qed\enddemo

\demo{Proof of  Proposition 8.1}
The implication (i) $\Rightarrow$ (ii) is trivial, since
$\pi_k(v)=\psi_v(p_u)$ is a particular character value for
$u=1^\infty21^{k-2}$. The statement (iii) follows from (ii) by
Lemma 8.4. In fact, we can split the initial sequence $\{w_n\}$
into two subsequences, $\{w'_n\}$ and $\{w''_n\}$, in such a way
that $d_1(w'_n)=2$ and $d_1(w''_n)\ne2$. Then we derive from
Lemma 8.4 that for both subsequences $\pi(w_n)\to0$, and (iii)
follows.

Now, (ii) follows from (iii) by Lemma 8.3, and (ii) implies (i)
by Lemma 8.2.
\qed\enddemo

\subhead General regularity conditions \endsubhead

We now find the  conditions for a sequence of linear functionals
on the spaces $R_n$ to converge to a functional on the limiting
space $R_\infty=\varinjlim R_n$.

\noindent
{\bf Definition \rm (Regularity of character sequences).}
Let $\psi_n$ be a linear functional on the graded component $R_n$
of the ring $R=\Bbb{Q}\langle X,Y\rangle$, for each
$n=1,2,\ldots$, and assume that the sequence converges pointwise
to a functional $\varphi$ on the ring $R$, in the sense that
$$
\lim_{n\to\infty} \psi_n(P\,X^{n-m}) = \varphi(P)
\tag 8.5
$$
for every $m\in\Bbb N$ and every polynomial $P\in R_m$. We call
such a sequence {\it regular}.

\smallskip
Our goal in this Section is to characterize the set of regular
sequences.

\smallskip\noindent
{\bf Definition \rm (Convergence of words).}
Let $\{w_n\}$ be a sequence of finite Fibonacci words, and assume
that the ranks $|w_n|$ tend to infinity as $n\to\infty$. We say
that $\{w_n\}$ {\it converges} to an infinite word $w$, iff the
$m$th letter $w_n(m)$ of $w_n$ coincides with the $m$th letter
$w(m)$ of the limiting word $w$ for almost all $n$ (i.e., for all
but finitely many $n$'s), and for all $m$.

Let us recall that an infinite word $w$ with $2$'s at positions
$d_1,d_2,\ldots$ is summable, if, and only if, the series
$\sum_{j=1}^\infty 1/d_j$ converges, or, equivalently, if the
product
$$
\pi(w) = \prod_{j:d_j\ge2} \Big(1 - {1 \over d_j}\Big) > 0
$$
converges.

Consider a sequence $w_1,w_2,\ldots$ of Fibonacci words
converging to a summable infinite word $w$. We denote by $w_n'$
the longest initial (rightmost) subword of $w_n$ identical with
the corresponding segment of $w$, and we call it {\it stable}
part of $w_n$. The remaining part of $w_n$ will be denoted by
$w_n''$, and referred to as {\it transient} part of $w_n$.

\smallskip\noindent
{\bf Definition \rm (Regularity Conditions).}
We say that a sequence of Fibonacci words $w_n\in \yf_n$ {\it
satisfies regularity conditions}, if either one of the following
two conditions holds:
$$
\lim_{n\to\infty} \pi(w_n) = 0;
\tag i
$$
or \smallskip \noindent
(ii)\; the sequence $w_n$ converges to a summable infinite word
$w$, and a strictly positive limit
$$
\beta = \pi^{-1}(w) \lim_{n\to\infty} \pi(w_n) > 0
\tag 8.6
$$
exists.

\proclaim{Theorem 8.5}
Assume that the regularity conditions hold for a sequence $w_n\in
\yf_n$. Then the character sequence $\psi_{w_n}$ is regular. If
the regularity condition (i) holds, then
$$
\lim_{n\to\infty} \psi_{w_n}(Q\,X^{n-m}) = \varphi_P(Q),
$$
and if regularity condition (ii) holds, then
$$
\lim_{n\to\infty} \psi_{w_n}(Q\,X^{n-m}) = \varphi_{\beta,w}(Q)
\tag 8.7
$$
for every polynomial $Q\in R_m$, $m=1,2,\ldots$.  Conversely,  if
the character sequence $\psi_{w_n}$ is regular, then the
regularity conditions hold for the sequence $w_n\in \yf_n$.
\endproclaim

This theorem will follow from Proposition 8.1 and the following
proposition:

\proclaim{Proposition 8.6}
Assume that a sequence $w_1,w_2,\ldots$ of Fibonacci words
converges to a summable infinite word $w$, and that there exists
a limit
$$
\beta = \pi^{-1}(w) \lim_{n\to\infty} \pi(w_n).
\tag 8.8
$$
Then
$$
\pi_k^{-1}(w) \lim_{n\to\infty} \pi_k(w_n) = \beta^k
\tag 8.9
$$
for every $k=2,3,\ldots$. More generally,
$$
\lim_{n\to\infty} \psi_{w_n}(p_u) = \beta^k\, \varphi_w(p_u)
\tag 8.10
$$
for every element $u\in 1^\infty \yf$ of essential rank $||u||=k$.
\endproclaim

\demo{Proof}
Let $m_n$ be the length of the stable part of the word $w_n$, and
note that $m_n\to\infty$. In the following ratio, the factors
corresponding to $2$'s in the stable part of $w_n$ cancel out,
$$
\frac{\pi(w_n)}{\pi(w)} =
\prod_{j:d_j(w_n)>m_n} \Big(1 - {1 \over d_j(w_n)}\Big)
\prod_{j:d_j(w)>m_n} \Big(1 - {1 \over d_j(w)}\Big)^{-1},
$$
and a similar formula holds for the functional $\pi_k$, $k=2,3,\ldots$:
$$
\frac{\pi_k(w_n)}{\pi_k(w)} =
\prod_{j:d_j(w_n)>m_n} \Big(1 - {k \over d_j(w_n)}\Big)
\prod_{j:d_j(w)>m_n} \Big(1 - {k \over d_j(w)}\Big)^{-1}.
\tag 8.11
$$
Consider the ratio
$$
\Big({\pi(w_n) \over \pi(w)}\Big)^k \Big/
\Big({\pi_k(w_n) \over \pi_k(w)}\Big) =
\prod_{j:d_j(w_n)>m_n} \frac
{\Big(1 - {1 \over d_j(w_n)}\Big)^k}
{\Big(1 - {k \over d_j(w_n)}\Big)}
\prod_{j:d_j(w)>m_n} \frac
{\Big(1 - {k \over d_j(w)}\Big)}
{\Big(1 - {1 \over d_j(w)}\Big)^k}.
$$
The second product in the right hand side is a tail of the
converging infinite product (since the word $w$ is summable),
hence converges to $1$, as $n\to\infty$. By (8.3), the first
product can also be estimated by a tail of a converging infinite
product,
$$
1 \le \prod_{j:d_j(w_n)>m_n} \frac
{\Big(1 - {1 \over d_j(w_n)}\Big)^k}
{\Big(1 - {k \over d_j(w_n)}\Big)} \le
\prod_{j:d_j>m_n} \Big(1 + {{k\choose2} \over (d_j-k)^2}\Big),
$$
hence converges to $1$, as well.

The proof of the formula (8.10) is only different in the way that
the ratio in the left hand side of (8.11) should be replaced by
$\psi_{w_n}(p_u)/\varphi_w(p_u)$.
\qed\enddemo

\demo{Proof of the Theorem 8.5}
It follows directly from Propositions 8.1 and 8.6 that the
regularity conditions for a sequence $w_1,w_2,\ldots$ imply
convergence of functionals $\psi_{w_n}$ to the Plancherel
character $\varphi_P$ in case (i), and to the character
$\varphi_{\beta,w}$ in case (ii). By the Corollary 7.5 we know
that the functions $\varphi_{\beta,w}$ are pairwise distinct, and
also different from the Plancherel functional $\varphi_P$.

Let us now assume that the sequence $\psi_{w_n}$ converges to a
limiting functional $\varphi$. We can choose a subsequence
$\psi_{w_{n_m}}$ in such a way that the corresponding sequence
$w_{m_n}$ converges digitwise to an infinite word $w$. If $w$ is
not summable, then $\varphi=\varphi_P$ coincides with the
Plancherel functional, and the part (i) of the regularity
conditions holds. Otherwise, we can also assume that the limit
(8.6) exists, and hence $\varphi=\varphi_{\beta,w}$ by
Proposition 8.6. Since the parameter $\beta$ and the word $w$ can
be restored, by Corollary 7.5, from the limiting functional
$\varphi$, the sequence $w_n$ cannot have subsequences converging
to different limits, nor can the sequence $\pi(w_n)$ have
subsequences converging to different limits. It follows, that the
regularity conditions are necessary. The Theorem is proved.
\qed\enddemo

In the following statement, $\Omega$ refers to the space defined
at the end of Section 3.

\proclaim{Theorem 8.7} The map
$$
(\beta, w) \mapsto \varphi_{\beta,w}, \quad P \mapsto \varphi_P
$$
is a homeomorphism of $\Omega$ onto the Martin boundary of $\yf$.
\endproclaim

\demo{Proof}
It follows from Corollary 7.5 that the map is an injection of
$\Omega$ into the Martin boundary, and from Theorem 8.5 that the
map is surjective. Furthermore, the proof that the map is a
homeomorphism is a straightforward variation of the proof of the
regularity statement, Theorem 8.5.
\qed\enddemo

\head \S 9. Concluding remarks \endhead

The Young-Fibonacci lattice, along with the Young lattice, are
the most interesting examples of differential posets. There is a
considerable similarity between the two graphs, as well as a few
severe distinctions.

Both lattices arise as Bratteli diagrams of increasing families
of finite dimensional semisimple matrix algebras, i.e., group
algebras of symmetric groups in case of Young lattice, and Okada
algebras in case of Young-Fibonacci graph. For every Bratteli
diagram, there is a problem of describing the traces of the
corresponding inductive limit algebra, which is well-known to be
intimately related to the Martin boundary construction for the
graph. The relevant fact is that indecomposable positive harmonic
functions, which are in one-to-one correspondence with the
indecomposable traces, form a part of the Martin boundary.

For the Young lattice the Martin boundary has been known for
several decades, and all of the harmonic functions in the
boundary are known to be indecomposable (extreme points). In this
paper we have found the Martin boundary for the Young-Fibonacci
lattice. Unfortunately, we still do not know which harmonic
functions in the boundary are decomposable (if any). The method
employed to prove indecomposability of the elements of the Martin
boundary of the Young lattice can not be applied to
Young-Fibonacci lattice, since the $K_0$-functor ring $R$ of the
limiting Okada algebra $\Cal F$ is not commutative, as it is in
case of the group algebra of the infinite symmetric group (in
this case it can be identified with the symmetric function ring).

Another natural problem related to Okada algebras is to find all
non-negative Markov traces. We plan to address this problem in
another paper.

\head Appendix \endhead

In this appendix, we survey a few properties of differential
posets introduced by R. Stanley in \cite{St1-3}. A more general
class of posets had been defined by S. Fomin in \cite{F1-3}. In
his terms differential posets are called self - dual graphs.

\subhead A.1 Definitions \endsubhead

A graded poset $\Gamma = \bigcup_{n=0}^\infty \Gamma_n$ is called
{\it branching diagram} (cf. \cite{KV}), if

(B1)\quad The set $\Gamma_n$ of elements of rank $n$ is finite
for all $n = 0,1,\ldots$

(B2)\quad There is a unique minimal element $\varnothing \in
\Gamma_0$

(B3)\quad There are no maximal elements in $\Gamma$.

One can consider a branching diagram as an extended phase space
of a non - stationary Markov chain, $\Gamma_n$ being the set of
admissible states at the moment $n$ and covering relations
indicating the possible transitions.

We denote the rank of a vertex $v \in \Gamma_n$ by $|v| = n$, and
the number of saturated chains in an interval $[u,v]\subset
\Gamma$ by $d(u,v)$.

Following \cite{St1}, we define an {\it $r$-differential poset}
as a branching diagram $\Gamma$ satisfying two conditions:

(D1)\quad If $u \ne v$ in $\Gamma$ then the number of elements
covered by $u$ and $v$ is the same as the number of elements
covering both $u$ and $v$

(D2)\quad If $v \in \Gamma$ covers exactly $k$ elements, then $v$
is covered by exactly $k+r$ elements of $\Gamma$.

Note that the number of elements in a differential poset covering
two distinct elements can be at most $1$. In this paper we focus
on $1$-differential posets.

For any branching diagram $\Gamma$ one can define two linear
operators in the vector space $\Fun(\Gamma)$ of  functions on
$\Gamma$ with coefficients in  $\Q$: the {\it creation operator}
$$
U(f)(v) = \sum_{w\colon v\nearrow w}f(w),
\tag A.1.1
$$
and the {\it annihilation operator}
$$
D(f)(v) = \sum_{u\colon u\nearrow v} f(u).
\tag A.1.2
$$
Identifying finitely supported functions on $\Gamma$ with formal
linear combinations of points of $\Gamma$ and vertices of
$\Gamma$ with the delta functions at the vertices, one can write
instead:
$$
U\, v = \sum_{w\colon v\nearrow w} w,
\tag A.1.1'
$$
$$
D\,v = \sum_{u\colon u\nearrow v} u.
\tag A.1.2'
$$

One can characterize $r$-differential posets  as branching
diagrams for which the operators $U,D$ satisfy the Weyl identity
$DU - UD = r\,I$.

\subhead A.2 Some properties of differential posets \endsubhead

We review below only a few identities we need in the main part of
the paper. For a general algebraic theory of differential posets
see \cite{St1-3}, \cite{Fo1-3}.  Assume here that $\Gamma$ is a
1-differential poset.

The first formula is well known:
$$
\sum_{w\colon v\nearrow w} d(\varnothing,w) = (n+1)\, d(\varnothing,v),
\quad v \in \Gamma.
\tag A.2.1
$$
\demo{Proof}
Let $d_n = \sum_{|v|=n} d(\varnothing,v)\, v \in \Fun(\Gamma)$. Then
$U\,d_n = d_{n+1}$ and (A.2.1) can be written as $D\,d_{n+1} =
(n+1)\, d_n$. This is trivial for $n=0$, and assuming $D\,d_n =
n\, d_{n-1}$ we obtain
$$
D\,d_{n+1} =
DU\,d_n = UD\,d_n + d_n = n\,U\,d_{n-1} + d_n = (n+1)\,d_n. \qed
$$
\enddemo

Our next result is a  generalization of (A.2.1).

\proclaim{Lemma A.2.2}
Let $\Gamma$ be a $1$-differential poset, and let $u \le v$ be
any vertices of ranks $|u|=k$, $|v|=n$. Then
$$
\sum_{w\colon v\nearrow w} d(u,w) -
\sum_{x\colon x\nearrow u} d(x,v) = (n-k+1)\,d(u,v).
\tag A.2.3
$$
\endproclaim

\demo{Proof}
Using the notation $d_n(u) = \sum_{|v|=n} d(u,v)\,v$, one can
easily see that $U\,d_n(u) = d_{n+1}(u)$ and that (A.2.3) can be
rewritten in the form
$$
D\,d_{n+1}(u) = \sum_{x\colon x\nearrow u} d_n(x) + (n-k+1)\,d_n(u).
$$
For $n=k-1$ the formula is true by the definition of $D$. By the
induction argument,
$$
\aligned
D\,d_{n+1}(u) & = DU\,d_n(u) = UD\,d_n(u) + d_n(u) = \\
& = U \left(\sum_{x\colon x\nearrow u} d_{n-1}(x) +
(n-k)\,d_{n-1}(u)\right) + d_n(u) = \\
& = \sum_{x\colon x\nearrow u} d_n(x) + (n-k+1)\, d_n(u).
\endaligned
$$
Note that (A.2.3) specializes to (A.2.1) in case $k=0$,
$u=\varnothing$.
\qed\enddemo

\subhead A.3 Plancherel transition probabilities on a
differential poset \endsubhead

It follows from (A.2.1) that the numbers
$$
p_{v,w} = {d(\varnothing,w) \over (n+1)\, d(\varnothing,v)};\quad
v\nearrow w,\; n=|v|,
\tag A.3.1
$$
can be considered as transition probabilities of a Markov chain
on $\Gamma$. Generalizing the terminology used in the particular
example of Young lattice (see \cite{KV}), we call (A.3.1) {\it
Plancherel transition probabilities}.

\proclaim{Lemma A.3.2}
Let $u \le v$ be vertices of ranks $|u|=k$, $|v|=n$ in a
$1$-differential poset $\Gamma$. Then the Plancherel probability
$p(u,v)$ to reach (by any path) the vertex $v$ starting with $u$
is
$$
p(u,v) = {k! \over n!}\, {d(u,v)\, d(\varnothing,v) \over
d(\varnothing,u)}.
\tag A.3.3
$$
\endproclaim

\demo{Proof}
We have to check that $\sum_v p(u,v)\, p_{v,w} = p(u,w)$. Since
$\sum_v d(u,v)\, d(v,w) = d(u,w)$, we obtain
$$
{k! \over n!} \sum_{|v|=n} {d(\varnothing,v) \over d(\varnothing,u)}
d(u,v)\,d(v,w){d(\varnothing,w) \over (n+1)\,d(\varnothing,v)} =
{k! \over (n+1)!}\, {d(u,w)\, d(\varnothing,w) \over
d(\varnothing,u)},
$$
and the Lemma follows.
\qed\enddemo



\Refs\widestnumber\no{KOO}

\ref\no D
\by J.~L.~Doob
\paper Discrete potential theory and boundaries
\jour J. of Math. and Mech.
\vol 8
\yr 1959
\pages 433-458
\endref

\ref \no E
\by\rm E.~Effros
\paper\it Dimensions and  $C^*$-algebras
\jour\rm CBMS Regional Conf. Series,
\publaddr\rm Providence
\vol 46
\yr 1981
\endref

\ref \no F1
\by \rm S.~V.~Fomin
\paper Generalized Robinson -- Schensted -- Knuth correspondence
\jour \it J. Soviet Math.
\vol \rm 41
\yr 1988
\pages 979-991
\endref

\ref \no F2
\by \rm S.~V.~Fomin
\paper Duality of graded graphs
\jour \it J. Soviet Math.
\vol \rm 41
\yr 1988
\pages 979-991
\endref

\ref \no F3
\by S.~V.~Fomin
\paper Schensted algorithms for dual graded graphs
\jour J. Algebraic Combinatorics
\vol 4
\yr 1995
\pages 5--45
\endref

\ref\no KV
\by S.~Kerov, A.~Vershik
\paper The Grothendieck Group of the Infinite Symmetric Group and
Symmetric Functions with the Elements of the $K_0$-functor theory
of AF-algebras
\inbook Representation of Lie groups and related topics
\bookinfo Adv. Stud. Contemp. Math.
\vol 7
\eds A.~M.~Vershik and D.~P.~Zhelobenko
\publ Gordon and Breach
\yr 1990
\pages 36--114
\endref

\ref\no KOO
\by S.~Kerov, A.~Okounkov and G.~Olshanski
\paper The Boundary of Young Graph with Jack Edge Multiplicities
\jour Kyoto University Preprint RIMS-1174
\yr 1997
\paperinfo q-alg/9703037
\endref

\ref\no M
\by I.~G.~Macdonald
\book Symmetric functions and Hall polynomials
\bookinfo 2nd edition
\publ Oxford University Press
\yr 1995
\endref

\ref \no Ok
\by \rm S.~Okada
\paper Algebras associated to the Young-Fibonacci lattice
\jour \it Trans. Amer. Math. Soc.
\vol \rm 346
\yr 1994
\pages 549--568
\endref

\ref\no Oku
\by A.~Yu.~Okounkov
\paper Thoma's theorem and representations of infinite bisymmetric
group
\jour Funct. Analysis and its Appl.
\vol 28
\yr 1994
\pages 101--107
\endref

\ref \no St1
\by \rm  R.~P.~Stanley
\paper Differential Posets
\jour \it J. Amer. Math. Soc.
\vol \rm 1
\yr 1988
\pages 919--961
\endref

\ref \no St2
\by \rm  R.~P.~Stanley
\paper Variations on Differential Posets
\inbook \it Invariant Theory and Tableaux, IMA Vol. Math. Appl.
\eds \rm D. Stanton
\yr 1988
\publ Springer Verlag
\publaddr New York
\pages 145--165
\endref

\ref \no St3
\by \rm  R.~P.~Stanley
\paper Further combinatorial properties of two Fibonacci lattices
\jour \it European J. Combin.
\vol \rm 11
\yr 1990
\pages 181--188
\endref

\ref\no T
\by\rm E.~Thoma
\paper\it Die unzerlegbaren, Positiv-definiten Klassen-funktionen der
abzahlbar un\-end\-li\-chen, sym\-me\-tri\-schen Gruppe
\jour\rm Math. Z.,
\vol 85
\yr 1964
\pages 40-61
\endref

\ref\no VK
\by A.~M.~Vershik, S.~V.~Kerov
\paper Asymptotic character theory of the symmetric group
\jour Funct. Analysis and its Appl.
\vol 15
\yr 1981
\pages 246--255
\endref

\endRefs
\enddocument
\end